# Szegö type results relative to a sequence of model spaces

BARUSSEAU BENOIT

I.M.B Université de Bordeaux I
351 cours de la Libération, 33405 Talence, France

*Email:* benoit.barusseau@math.u-bordeaux1.fr

**Abstract**

The classical Szegö theorem can be stated in terms of the sequence of model spaces $(K_{z^N})_{N \in \mathbb{N}}$. In this article, we are interested in the generalization of the Szegö theorem in the case of the sequence $(K_{B^N})_{N \in \mathbb{N}}$ where $B$ is a finite Blaschke product.

**A.M.S. subject classification:** Operator theory

The Hardy space of the unit circle $\mathbb{T}$ is defined to be $H^2 = \{f \in L^2(\mathbb{T}, dv), \hat{f}(n) = 0 \text{ if } n < 0\}$ where $\hat{f}(n)$ denotes the $n$-th order Fourier coefficient of $f$. $H^2$ is a Hilbert space with scalar product :

$$<f,g> = \int_{\mathbb{T}} f(\xi) \overline{g}(\xi) \frac{d\xi}{2\pi} \quad \forall f, g \in H^2.$$

We recall that $H^2$ has the reproducing kernel :

$$K_\lambda(z) = \frac{1}{1 - \overline{\lambda} z}.$$

By this we mean that for all $\lambda \in \mathbb{D}$, $K_\lambda$ is an element of $H^2$ and $f(\lambda) = <f, K_\lambda>$, for all $f \in H^2$. We denote $k_z = \frac{K_z}{\|K_z\|}$ the normalized reproducing kernel such that

$$k_z(\lambda) = \frac{(1 - |z|^2)^{1/2}}{1 - \overline{\lambda} z}.$$

Let $\lambda \in \mathbb{D}$ and $b_\lambda(\xi) = \frac{\xi - \lambda}{1 - \overline{\lambda} \xi}$. The Blaschke product $B$ with zero sequence

$$Z(B) = \{\lambda_1, ..., \lambda_n\}$$

is defined to be

$$B(z) = \prod_{i=1}^n b_{\lambda_i}(z) = \prod_{i=1}^n \frac{z - \lambda_i}{1 - \overline{\lambda_i} z}.$$

In particular $B$ is an inner function, that is $|B| = 1$ almost everywhere on $\mathbb{T}$.

The Toeplitz operator with symbol $\varphi \in L^2$ is defined on the dense subset $H^\infty$ of $H^2$ to be $T_\varphi(g) = P^{H^2}(\varphi g)$. It is well known that $T_\varphi$ can be extended to a bounded operator on $H^2$ if and only if $\varphi \in L^\infty$ (see [3])

Let $\varphi \in L^\infty$ and let $\theta \in H^\infty = L^\infty(\mathbb{T}) \cap H^2$ be an inner function, thus $\theta H^2$ is a closed subspace of $H^2$ and we denote by $K_\theta$ the model space $H^2 \ominus \theta H^2 = H^2 \cap (\theta H^2)^\perp$. It is well know that the spaces $\theta H^2$ are shift invariant and any shift invariant subspace of $H^2$ is of that form. Thus it is easy to see that any model space is backward shift invariant. We denote by $P_{K_\theta}$ is the orthogonal projection from $H^2$ to $K_\theta$ and by $T_{\varphi, K_\theta}$ denotes the compression to $K_\theta$ of the classical Toeplitz operator $T_\varphi$, this is $T_{\varphi, K_\theta} = P_{K_\theta} T_\varphi P_{K_\theta}$. Thus the matrix of $T_{\varphi, K_{z^n}}$ in the basis $\{1, z, ..., z^{n-1}\}$ is a classical finite dimensional Toeplitz matrix of size $n$.

If $\theta$ is a finite Blaschke product then $K_{\theta^n}$ is a finite dimensional space and for any $\varphi \in L^2$, $T_{\varphi, K_\theta}$ is a continuous operator.

Recently D. Sarason [6] discussed the properties (finite rank, compacity...) of such compressions of Toeplitz operators to some model spaces. We are interested in establishing a type of Szegö theorem for these operators.

In the first section, we recall the classical Szegö theorem. Next we study the space $K_{B^n}$ and give a basis which will be useful in calculating the matrix for a certain function class $M_1 \subset L^2(\mathbb{T})$. Finally, considering this space, we will be able to give a Szegö-type result for compressions of Toeplitz operator to model spaces.

# 1 The classical Szegö theorem.

**Definition 1.**

1. Let $(A_n)_{n \in \mathbb{N}}$ be a sequence of matrices of size $n \times n$ and $f$ a measurable function on $\mathbb{T}$. $A_n$ is said to be distributed as $f$ in the sense of the eigenvalues if and only if for any continuous function $G$ with compact support we have

$$\lim_{n \to \infty} \frac{1}{n} \sum_{k=1}^n G(\lambda_k(A_n)) = \frac{1}{2\pi} \int_\mathbb{T} G(f(\xi)) d\xi$$

   where $(\lambda_k(A_n))_{1 \leqslant k \leqslant n}$ denotes the sequence of the eigenvalues of $A_n$.

2. $A_n$ is said to be distributed as $f$ in the sense of the singular values if and only if for any continous function $G$ with compact support, we have

$$\lim_{n \to \infty} \frac{1}{n} \sum_{k=1}^n G(\sigma_k(A_n)) = \frac{1}{2\pi} \int_\mathbb{T} G(f(\xi)) d\xi$$

   where $(\sigma_k(A_n))_{1 \leqslant k \leqslant n}$ denotes the sequence of the singular values of $A_n$.

We shall use the $L^2$ version [8] of the classical Szegö theorem and the theorem of Avram [1] and Parter [5] concerning the singular values. We can find demonstrations in the continuous case using the operator theory and the $C^*$-algebras in [2]. Let us recall this fundamental theorem :

**Theorem 2.** *[8] Let $f \in L^2(\mathbb{T})$, the sequence of Toeplitz truncated matrices of size n is distributed as $|f|$ in the sense of the singular values.. If $f$ is real valued then this sequence is distributed as $f$ in the sense of the eigenvalues.*

## 2 General sequence of model spaces $K_{B^n}$.

### 2.1 Vocabulary

We recall a fact and introduce vocabulary. It is well known (see [3] or [4]) that $K_\theta$ is finite dimensional if and only if there exists $B$ a finite Blaschke product such that $\theta = cB$ where $c \in \mathbb{T}$. In order to consider the eigenvalues of the compression of Toeplitz operators to model spaces, we only consider finite dimensional model spaces. Thus we only consider the case of finite Blaschke products.

**Definition 3.** *Let $B \in H^2$ be a finite Blaschke product.*

1. *$T_f$ is said to be distributed as $g$ in the sense of the eigenvalues relative to the sequence of model spaces $(K_{B^n})_{n \in \mathbb{N}}$ if and only if, for all continuous functions $G$ with compact support, we have*

$$\lim_{n \to \infty} \frac{1}{\dim K_{B^n}} \sum_{k=1}^{\dim K_{B^n}} G(\lambda_k(T_{f,K_{B^n}})) = \frac{1}{2\pi} \int_{\mathbb{T}} G(g(\xi)) d\xi$$

2. *$T_f$ is said to be distributed as $g$ in the sense of the singular values relative to the sequence of model spaces $(K_{B^n})_{n \in \mathbb{N}}$ if and only if, for all continuous functions $G$ with compact support, we have*

$$\lim_{n \to \infty} \frac{1}{\dim K_{B^n}} \sum_{k=1}^{\dim K_{B^n}} G(\lambda_k(T_{f,K_{B^n}})) = \frac{1}{2\pi} \int_{\mathbb{T}} G(g(\xi)) d\xi$$

Now, we can state theorem 2 in the following way : If $f \in L^2$ then $T_f$ is distributed as $|f|$ in the sense of the singular values relative to the sequence of model spaces $(K_{z^N})_{N \in \mathbb{N}}$. If $f$ is real valued then $T_f$ is distributed as $f$ in the sense of the eigenvalues relative to the sequence of model space $(K_{z^N})_{N \in \mathbb{N}}$.

Now, we study the case of a general finite Blaschke product.

### 2.2 Basic properties of $K_{B^N}$

While many bases have been proposed for model spaces, the Malmquist basis [4] is particularly helpful when dealing with our problem. We only need to choose an ordering of the basis elements appropriate to our problem.

Let $B$ be a finite Blaschke product. We denote by $m(B, \lambda_i)$ the multiplicity of $\lambda_i$ as a zero of $B$. Thus there exists $p \in \mathbb{N}$ and distinct $\lambda_i$'s such that

$$B = \prod_{i=1}^{p} b_{\lambda_i}^{m(B,\lambda_i)}$$

For any $1 \leqslant j, k \leqslant p$, $0 \leqslant r \leqslant m(B, \lambda_j) - 1$ and $0 \leqslant s \leqslant m(B, \lambda_k) - 1$, we denote

$$e_j^r = b_{\lambda_j}^r \prod_{i=1}^{j-1} b_{\lambda_i}^{m(B,\lambda_i)} k_{\lambda_j}.$$

and consider the basis

$$S = \left( (e_j^r)_{r=0\ldots m(B,\lambda_j)-1} \right)_{j=1\ldots p}.$$

Thus we remark that denoting $S_j = \{k_{\lambda_j}, b_{\lambda_j}k_{\lambda_j}, ...b_{\lambda_j}^{m(B,\lambda_j)-1}k_{\lambda_j}\}$ for $j \in \{1, ..., p\}$, we have

$$S = S_1 \bigcup \left(b_{\lambda_1}^{m(B,\lambda_1)}S_2\right) \bigcup \left(b_{\lambda_1}^{m(B,\lambda_1)}b_{\lambda_2}^{m(B,\lambda_2)}S_3\right) \bigcup ... \bigcup \left(\prod_{i=1}^{p-1} b_{\lambda_i}^{m(B,\lambda_i)} S_p\right)$$

where $fA$ is the set $\{fg, g \in A\}$.

In the following we denote $S$ vertically: $\left\{\begin{array}{c} e_1^0 \\ \vdots \\ e_p^{m(B,\lambda_p)-1} \end{array}\right\}$.

Now, we show that $S$ is an orthonormal basis for $K_B$ and generalize to get an orthonormal one for $K_{B^n}$.

**Lemma 4.** *For all $1 \leqslant j, k \leqslant p$, $0 \leqslant r \leqslant m(B, \lambda_j) - 1$ and $0 \leqslant s \leqslant m(B, \lambda_k) - 1$, we have*

$$<e_j^r, e_k^s> = \delta_{j,k}\delta_{r,s}.$$

**Proof.** We remark that $e_j^r = b_{\lambda_j}^r \prod_{i=1}^{j-1} b_{\lambda_i}^{m(B,\lambda_i)} k_{\lambda_j}$. Without loss of generality, we can suppose $j \geqslant k$. If $j > k$

$$<b_{\lambda_j}^r \prod_{i=1}^{j-1} b_{\lambda_i}^{m(B,\lambda_i)} k_{\lambda_j}, b_{\lambda_k}^s \prod_{i=1}^{k-1} b_{\lambda_i}^{m(B,\lambda_i)} k_{\lambda_k}> = <b_{\lambda_j}^r \prod_{i=k}^{j-1} b_{\lambda_i}^{m(B,\lambda_i)} k_{\lambda_j}, b_{\lambda_k}^s k_{\lambda_k}>$$

$$= <b_{\lambda_j}^r b_{\lambda_k}^{m(B,\lambda_k)-s} \prod_{i=k+1}^{j-1} b_{\lambda_i}^{m(B,\lambda_i)} k_{\lambda_j}, k_{\lambda_k}>$$

where we set $\prod_{i=k+1}^{j-1} b_{\lambda_i}^{m(B,\lambda_i)} = 1$ for $k+1 > j-1$. Since $1 \leqslant m(B, \lambda_k) - s$, the function on the left hand side has a zero in $\lambda_k$ and so

$$<b_{\lambda_j}^r b_{\lambda_k}^{m(B,\lambda_k)-s} \prod_{i=k+1}^{j-1} b_{\lambda_i}^{m(B,\lambda_i)} k_{\lambda_j}, k_{\lambda_k}> = 0.$$

Suppose $k = j$, then

$$<b_{\lambda_j}^r \prod_{i=1}^{j-1} b_{\lambda_i}^{m(B,\lambda_i)} k_{\lambda_j}, b_{\lambda_k}^s \prod_{i=1}^{k-1} b_{\lambda_i}^{m(B,\lambda_i)} k_{\lambda_k}> = <b_{\lambda_j}^r k_{\lambda_k}, b_{\lambda_k}^s k_{\lambda_k}> = \delta_{s,r}.$$

□

Now since the space $K_B$ has dimension $|Z(B)|$, $S$ is an orthonormal basis for $K_B$. So denoting

$$\Sigma = \{S, BS, ..., B^{n-1}S\} = \left(\left(\left(B^k e_j^r\right)_{k=0...n-1}\right)_{r=0...m(B,\lambda_j)-1}\right)_{j=1...p},$$

it is easy to see that $\Sigma$ is a basis for $K_{B^n}$.

**Proposition 5.** *[4] (Malmquist Basis) For all $n \in \mathbb{N}$, $\Sigma$ is an orthonormal basis of $K_{B^n}$.*

**Proof.** It follows from the fact that $B^k S$ and $B^p S$ are orthogonal for $k \neq p$. Indeed if $f \in B^k S$ and $g \in B^p S$ and $k > p$ then

$$<fB^k, gB^p> = <fB^{k-p}, g>$$

and since $fB^{k-p} \in BH^2$ and $g \in S \subset K_B$ we have the result. □

We have the following corollary (the proof is obvious).

**Corollary 6.** *The sequence $(K_{B^n})_{n \in \mathbb{N}}$ verifies*

1. *$(K_{B^n})_{n \in \mathbb{N}}$ is an increasing sequence (with respect to the order $\subset$ in $H^2$).*
2. *$(\dim K_{B^n})_n$ is strictly increasing and $\dim K_{B^{n+1}} - \dim K_{B^n} = |Z(B)|$.*

We can write $\Sigma$ as

$$\Sigma = \left\{ \begin{array}{c} e_1^0, Be_1^0, ..., B^{n-1}e_1^0, \\ \vdots \\ e_1^{m(B,\lambda_1)-1}, Be_1^{m(B,\lambda_1)-1}..., B^{n-1}e_1^{m(B,\lambda_1)-1}, \\ e_2^0, Be_2^0, ..., B^{n-1}e_2^0, \\ \vdots \\ e_p^0, Be_p^0..., B^{n-1}e_p^0, \\ \vdots \\ e_p^{m(B,\lambda_p)-1}, Be_p^{m(B,\lambda_p)-1}..., B^{n-1}e_p^{m(B,\lambda_p)-1} \end{array} \right\}.$$

We denote by $\Sigma_i$ the $i$-th line of $\Sigma$. Thus if we consider elements of the form $B^i e_j^r$ and $B^l e_k^s$ then they are in distincts $\Sigma_l$'s if and only if $j \neq k$ or $r \neq s$. In the following, we are interested in the matrix of $T_{\varphi,B^n}$ relative to the blocs $\Sigma_l$. $\Sigma$ is now identified with the vector $(\Sigma_1, ..., \Sigma_{n \times |Z(B)|})^t$. Thus $\Sigma$ is a basis of $K_{B^n}$ and it is now listed in the order

$$\left( \left( \left( B^k e_j^r \right)_{j=1...p} \right)_{k=0...n-1} \right)_{r=0...m(B,\lambda_j)-1}.$$

## 2.3 Szegö type result.

We need to do some calculations

**Lemma 7.** *Let $B$ be a finite Blaschke product with $B = \prod_{q=0}^p b_{\lambda_q}^{m(B,\lambda_q)}$.*

*For all $i, l \in \mathbb{N}$, $n \in \mathbb{Z}$, $1 \leqslant j, k \leqslant p$, $0 \leqslant r \leqslant m(B,\lambda_j)-1$ and $0 \leqslant s \leqslant m(B,\lambda_k)-1$, we have*

$$< B^n B^i e_j^r, B^l e_k^s > = \delta_{n+i,l} \delta_{j,k} \delta_{r,s}$$

**Proof.**
If $n \geqslant -i$ then $B^n B^i e_j^r \in B^{n+i} S$ thus $< B^n B^i e_j^r, B^l e_k^s > = \delta_{n+i,l} \delta_{j,k} \delta_{r,s}$.
If $n < -i$ then $< B^n B^i e_j^r, B^l e_k^s > = < B^i e_{j,r}, B^{l-n} e_{k,s} >$ where $n-l \geqslant 0$ and by the same reasoning, we have $< B^i e_j^r, B^{l-n} e_k^s > = \delta_{i,l-n} \delta_{j,k} \delta_{r,s} = \delta_{n+i,l} \delta_{j,k} \delta_{r,s}$. $\square$

Now we can give the matrix of $T_{\phi, K_{B^n}}$ relative to $\Sigma$. We begin by considering symbols in $M_1 = \overline{\text{Vect}(B^t, t \in \mathbb{Z})}^{L^2}$.

**Theorem 8.** *Let $\phi = \sum_{t \in \mathbb{Z}} a_t B^t \in M_1$ and $\varphi(z) = \sum_{|t| \leqslant n} a_t z^t$. The matrix of $T_{\phi, K_{B^n}}$ relative to $\Sigma$ (still denoted by $T_{\phi, K_{B^n}}$) is a diagonal block Toeplitz matrix of size $|Z(B)| \times n$. More precisely*

$$T_{\phi, K_{B^n}} = \begin{pmatrix} T_{\varphi, K_{z^n}} & & 0 \\ & \ddots & \\ 0 & & T_{\varphi, K_{z^n}} \end{pmatrix}$$

**Proof.** Let $\phi = \sum_{t \in \mathbb{Z}} a_t B^t$. We have to calculate for any $i, l \leqslant n-1$

$$<T_{\phi, K_{B^n}} B^i e_j^r, B^l e_k^s> = \sum_{t \in \mathbb{Z}} a_t <B^t B^i e_j^r, B^l e_k^s>$$

where $1 \leqslant j, k \leqslant p$, $0 \leqslant r \leqslant m(B, \lambda_j) - 1$ and $0 \leqslant s \leqslant m(B, \lambda_k) - 1$. Thus, using lemma 7, we have

$$\sum_{t \in \mathbb{Z}} a_t <B^t B^i e_j^r, B^l e_k^s> = a_{l-i} \delta_{j,k} \delta_{r,s}. \tag{1}$$

We denote by $A_{u,v}$ the blocs corresponding to the subspaces $\Sigma_u$ and $\Sigma_v$, that is

$$\begin{array}{c} & \Sigma_1 & \cdots & \Sigma_{|Z(B)|} \\ \Sigma_1 \\ \vdots \\ \Sigma_{|Z(B)|} \end{array} \begin{pmatrix} A_{1,1} & \cdots & A_{1,|Z(B)|} \\ \vdots & \ddots & \vdots \\ A_{|Z(B)|,1} & \cdots & A_{|Z(B)|,|Z(B)|} \end{pmatrix}$$

The term $\delta_{j,k}\delta_{r,s}$ in (1) means that $A_{u,v}$ is the zero bloc if $u \neq v$. If $j = k$ and $r = s$ then we are on a diagonal bloc and the term $a_{l-i}$ means that this bloc is Toeplitz of size $n$. Thus, in the basis $\Sigma$, we have the following matrix

$$\begin{array}{c} & \Sigma_1 & \cdots & \Sigma_{|Z(B)|} \\ \Sigma_1 \\ \vdots \\ \Sigma_{|Z(B)|} \end{array} \begin{pmatrix} T_{\varphi, K_{z^n}} & 0 & 0 \\ 0 & \ddots & 0 \\ 0 & 0 & T_{\varphi, K_{z^n}} \end{pmatrix}$$

$\square$

In the next section, we show that the space $M_1$ consists of all functions $\phi = \sum_{t \in \mathbb{Z}} a_t B^t$ where $(a_t)_{t \in \mathbb{Z}} \in \ell^2(\mathbb{Z})$. Then in the final section, we will be able to state a Szegö type result.

## 3 The space $M_1$

### 3.1 The general case.

**Lemma 9.** *Let $h \colon \mathbb{T} \to [0, \infty]$ be measurable on $\mathbb{T}$ and $B$ a finite Blaschke product. We have the following inequality*

$$\int_{\mathbb{T}} h(\xi) d\xi \leqslant \int_{\mathbb{T}} h(B(\xi)) \times |B'(\xi)| d\xi \tag{2}$$

**Proof.** In order to use the change of variable $\xi = B(e^{i\theta})$, we denote by $j$ the function $\theta \longmapsto B(e^{i\theta})$ and by $l$ the principal argument of $j$. Denoting by Log the principal value of the logarithm and Arg the principal argument, one has $\text{Log}(z) = \ln(|z|) + i \, \text{Arg}(z)$ and $\text{Log}(e^{i\theta}) = i\theta$ where $\theta \in \,]-\pi, \pi]$. Now, since $B$ has modulus 1, $l$ verifies the following equation

$$\frac{1}{i} \text{Log}(B(e^{i\theta})) = l(\theta)$$

for all $\theta \in \,]-\pi, \pi]$. If $B$ has one zero $B = b_\lambda$ then

$$l'(\theta) = \frac{i \, e^{i\theta} b_\lambda'(e^{i\theta})}{i \, b_\lambda(e^{i\theta})} = e^{i\theta} \frac{(1-|\lambda|^2)}{(1-\bar{\lambda} e^{i\theta})(e^{i\theta}-\lambda)} = \frac{(1-|\lambda|^2)}{|1-\lambda e^{i\theta}|^2}.$$

Thus considering the given finite Blaschke product $B$, we have

$$l'(\theta) = \sum_j \frac{(1-|\lambda_j|^2)}{|1-\lambda_j e^{i\theta}|^2} > 0.$$

Thus since $j$ is surjective, continuously differentiable on $]-\pi,\pi]$ and since its argument has a strictly increasing derivative, there exists $I \subset ]-\pi,\pi]$ such that $j$ is a continuously differentiable bijection from $I$ to $\mathbb{T}$. Thus considering the change of variable $\xi = j(\theta)$ we obtain

$$\int_{\mathbb{T}} h(\xi)d\xi = \int_{j(I)} h(\xi)d\xi = \int_I h(j(\theta)) \times |j'(\theta)|d\theta$$

moreover $j'(\theta) = i\, e^{i\theta} B'(e^{i\theta})$, which implies

$$\int_I h(j(\theta)) \times |j'(\theta)|d\theta = \int_I h(B(e^{i\theta})) \times |B'(e^{i\theta})|d\theta.$$

It is clear that $e^{iI} \subset \mathbb{T}$ so by positivity of the integrated function, we have

$$\int_I h(B(e^{i\theta})) \times |B'(e^{i\theta})|d\theta \leqslant \int_{\mathbb{T}} h(B(\xi)) \times |B'(\xi)|d\xi.$$

$\square$

**Remark 10.** This proof uses only the fact that the winding number of a Blaschke product is the number of its zeros.

**Proposition 11.**

1. If $(a_n)_{n\in\mathbb{Z}} \in \ell^2(\mathbb{Z})$ then $\sum_{n\in\mathbb{Z}} a_n B^n \in M_1$.

2. Let us define $\Gamma: \ell^2 \to M_1 \subset L^2$ by

$$\Gamma(a_n) = \sum_{n\in\mathbb{Z}} a_n B^n.$$

$\Gamma$ is a bicontinuous isomorphism. Precisely, we have

$$\frac{1 - |\prod_{\lambda \in Z(B)} \lambda|}{2} \|\Gamma(a_n)\|^2 \leqslant \|(a_n)\|_2^2 \leqslant \sup_{\xi \in \mathbb{T}} |B'(\xi)| \times \|\Gamma(a_n)\|^2.$$

**Proof.** 1) We see that

$$<\sum_{n\in\mathbb{Z}} a_n B^n, \sum_{n\in\mathbb{Z}} a_n B^n> = \sum_{n,k\in\mathbb{Z}} a_n a_k <B^n, B^k>,$$

Now, denoting $\delta = B(0)$, we have $|\delta| = \left|\prod_{\lambda \in Z(B)} \lambda\right| < 1$. Using the fact that 1 is the reproducing kernel associated with 0, we have $<B^n, B^k> = <B^{n-k}, 1> = B^{n-k}(0) = \delta^{n-k}$ for $n \geqslant k$ thus

$$<\sum_{n\in\mathbb{Z}} a_n B^n, \sum_{k\in\mathbb{Z}} a_k B^k> = \sum_{n\geqslant k} a_n a_k \delta^{n-k} + \sum_{k>n} a_n a_k \bar{\delta}^{n-k}. \qquad (3)$$

We can write

$$\sum_{n\geqslant k} a_n a_k \delta^{n-k} = \sum_{u\geqslant 0} \left(\sum_{n\in\mathbb{Z}} a_n a_{n-u}\right) \delta^u. \qquad (4)$$

Moreover for all $u \in \mathbb{Z}$

$$\left| \sum_{n \in \mathbb{Z}} a_n a_{n-u} \right| = \left| < \sum_{i \in \mathbb{Z}} a_i z^i, z^u \sum_{i \in \mathbb{Z}} a_i z^i > \right| \leqslant \|z^u\| \sum_{i \in \mathbb{Z}} |a_i|^2 = \sum_{i \in \mathbb{Z}} |a_i|^2. \quad (5)$$

Thus by equations (4) and (5), we have

$$\left| \sum_{n \geqslant k} a_n a_k \delta^{n-k} \right| \leqslant \sum_{u \geqslant 0} \|(a_n)\|_2^2 . |\delta|^u = \|(a_n)\|_2^2 \frac{1}{1-|\delta|}.$$

By the same reasoning on $\sum_{k>n} a_n a_k \bar{\delta}^{n-k}$, we have

$$\left| \sum_{k > n} a_n a_k \bar{\delta}^{n-k} \right| \leqslant \sum_{u > 0} \|(a_n)\|_2^2 . |\delta|^u \leqslant \|(a_n)\|_2^2 \frac{1}{1-|\delta|} \quad (6)$$

Thus combining (3), (5) and (6), we get

$$\|\Gamma(a_n)\|_2^2 \leqslant \|(a_n)\|_2^2 \frac{2}{1-|\delta|}.$$

2) Since $\{B^n, n \in \mathbb{Z}\}$ is a basis of $M_1$, $\Gamma$ is obviously an isomorphism. 1) shows that $\Gamma$ is continuous. Suppose $\Gamma(a_n) \in M_1 \subset L^2$. By lemma 9, we have

$$\|(a_n)\|_2^2 = \int_{\mathbb{T}} \left| \sum_{n \in \mathbb{Z}} a_n \xi^n \right|^2 d\xi \leqslant \int_{\mathbb{T}} |\varphi(B(\xi))|^2 \times |B'(\xi)| d\xi \leqslant \sup_{\xi \in \mathbb{T}} |B'(\xi)| \times \|\Gamma(a_n)\|_2^2.$$

Moreover $\sup_{\xi \in \mathbb{T}} |B'(\xi)| < \infty$ and since $B$ is a finite Blaschke product, $B$ is holomorphic in a neighborhood of $\overline{\mathbb{D}}$, thus $(a_n) \in \ell^2(\mathbb{Z})$. $\square$

In the case where $|Z(B)| = 1$, we show that $M_1$ is equal to $L^2$. Moreover in the next section, we will see that our Szegö type result can be expressed in a simple way.

**Proposition 12.** *If $\lambda \in \mathbb{D}$ and $B = b_\lambda$ then $M_1 = L^2$.*

**Proof.** We define the operator $\Gamma_\lambda = \Gamma \colon L^2 \longrightarrow L^2$ by $\Gamma_\lambda(\sum_{n \in \mathbb{Z}} a_n) = \sum_{n \in \mathbb{Z}} a_n b_\lambda^n$ then proposition 11 gives

$$\|\Gamma_\lambda\|^2 \leqslant \frac{2}{1-|\lambda|}.$$

It is easy to see that $b_\lambda \circ b_{-\lambda}(z) = b_{-\lambda} \circ b_\lambda(z) = z$ thus $\Gamma_\lambda \Gamma_{-\lambda} = \Gamma_{-\lambda} \Gamma_\lambda = I$ and $\|\Gamma_{-\lambda}\|^2 \leqslant \frac{2}{1-|\lambda|}$ thus $M_1 = \Gamma_\lambda(L^2) = L^2$. $\square$

## 4 The results

Proposition 11 shows that $\phi = \sum_{n \in \mathbb{Z}} a_n B^n \in M_1$ if and only if $\varphi(z) = \sum_{n \in \mathbb{Z}} a_n z^n \in L^2(\mathbb{T})$. Finally, as a corollary of theorem 8, and using the classical Szegö theorem, we obtain the following theorem.

**Theorem 13.** *Let $\phi = \sum_{n \in \mathbb{Z}} a_n B^n \in M_1$, denoting $\varphi(z) = \sum_{n \in \mathbb{Z}} a_n z^n \in L^2(\mathbb{T})$ then $T_\phi$ is distributed as $|\varphi|$ in the sense of the singular values relative to the sequence of model spaces $(K_{B^n})_{n \in \mathbb{N}}$. If $\phi$ is real valued then $T_\phi$ is distributed as $\varphi$ in the sense of the eigenvalues relative to the sequence of model spaces $(K_{B^n})_{n \in \mathbb{N}}$.*

**Proof.** By theorem 8, the matrix of $T_{\phi,K_{B^n}}$ relative to the basis $\Sigma$ is

$$T_{\phi,K_{B^n}} = \begin{pmatrix} T_{\varphi,K_{z^n}} & 0 & 0 \\ 0 & \ddots & 0 \\ 0 & 0 & T_{\varphi,K_{z^n}} \end{pmatrix}$$

where this matrix has $|Z(B)|$ blocs and the dimension of $K_{B^n}$ is $n \times |Z(B)|$. Thus the singular values of $A_n$ are just those of $T_{\varphi,K_{z^n}}$ with multiplicity $|Z(B)|$. So for all continuous functions $G$ with compact support, we can write

$$\frac{1}{\dim K_{B^n}} \sum_{k=1}^{\dim K_{B^n}} G(\sigma_k(T_{\phi,K_{B^n}})) = \frac{|Z(B)|}{n \times |Z(B)|} \sum_{k=1}^{n} G(\sigma_k(T_{\varphi,K_{z^n}}))$$
$$= \frac{1}{n} \sum_{k=1}^{n} G(\sigma_k(T_{\varphi,K_{z^n}}))$$

and theorem 2 allows us to conclude. The same reasoning gives the distribution of the eigenvalues in the real case, it suffices to remark that if $\phi$ is real valued so is $\varphi$. □

In fact we could use results on block Toeplitz matrices (see [2] and [7]) and we will discuss that point at the end of this article.

As an important corollary, we can state a simple form of theorem 13 in the case $|Z(B)| = 1$.

**Theorem 14.** *Let $f \in L^2$ be a real valued function. Then $T_f$ is distributed as $f \circ b_{-\lambda}$ in the sense of the eigenvalues relative to the sequence of model spaces $(K_{b_\lambda^n})_{n \in \mathbb{N}}$, so, for all continuous functions $G$ with compact support, we have*

$$\lim_{n \to \infty} \frac{1}{\dim K_{b_\lambda^n}} \sum_{k=1}^{\dim K_{b_\lambda^n}} G(\lambda_k(T_{f,K_{b_\lambda^n}})) = \frac{1}{2\pi} \int_{\mathbb{T}} G((f \circ b_{-\lambda}))(\xi) d\xi$$

*and $T_f$ is distributed as $|f \circ b_{-\lambda}|$, in the sense of the singular values relative to the sequence of model spaces $(K_{b_\lambda^n})_{n \in \mathbb{N}}$. Thus, if $G$ is a continuous function with compact support, we have*

$$\lim_{n \to \infty} \frac{1}{\dim K_{b_\lambda^n}} \sum_{k=1}^{\dim K_{b_\lambda^n}} G(\sigma_k(T_{f,K_{b_\lambda^n}})) = \frac{1}{2\pi} \int_{\mathbb{T}} G(|f \circ b_{-\lambda}|)(\xi) d\xi$$

The case of general symbols is still open. To solve it we would like to express symbol $\phi \in L^2$ so that the matrix of $T_{\phi,K_{B^n}}$ in the basis $\Sigma$ allows us to apply the $L^2$ version of the Szegö theorem.

We give some simple remarks exploring the case $|Z(B)| = 2$ and $B = b_{\lambda_1} b_{\lambda_2}$. In fact we need to complete the set $\{B^n, n \in \mathbb{Z}\}$ to obtain a basis for $L^2(\mathbb{T})$. We first observe that the set

$$\{k_{\lambda_1}, b_{\lambda_1} k_{\lambda_2}, B k_{\lambda_1}, B b_{\lambda_1} k_{\lambda_2}, B^2 k_{\lambda_1}, B^2 b_{\lambda_1} k_{\lambda_2} ...\}$$

is an orthonormal basis for $H^2$. (Indeed if $g \in H^2$ is orthogonal to any element of the previous set then $<g, k_{\lambda_1}> = 0$ so $g = b_{\lambda_1} g_2$ where $g_2 \in H^2$ and $<g, b_{\lambda_1} k_{l_2}> = <g_2, k_{\lambda_2}> = 0$ implies that $g_2 = b_{\lambda_2 k_{\lambda_2}}$ and $g = Bg_3$. Using the same reasoning, we obtain that $\lambda_1$ and $\lambda_2$ are zeros for $g_3$ and, by induction, $\lambda_1$ is a zero for $g$ with infinite order so $g = 0$.)

Now considering rational functions, it is clear that

$$k_{\lambda_1}, b_{\lambda_1}k_{\lambda_2} \in \mathrm{Span}(1, B, b_{\lambda_1}) = \left\{\frac{P}{(1-\overline{\lambda_1}X)(1-\overline{\lambda_2}X)}, P \in \mathbb{C}_2[X]\right\}.$$

Thus using the multiplication by $B$ we see that $\{B^n b_{\lambda_1}, n \in \mathbb{N}\} \bigcup \{B^n, n \in \mathbb{Z}\}$ is a basis of $H^2(\mathbb{T})$. Now using the equation $L^2(\mathbb{T}) = H^2(\mathbb{T}) \oplus \overline{H_0^2(\mathbb{T})}$ where $H_0^2(\mathbb{T}) = \{f \in H^2(\mathbb{T}), f(0) = 0\}$, we can write any element $\phi$ of $L^2(\mathbb{T})$ as

$$\phi = \sum_{n \in \mathbb{Z}} a_n B^n + \sum_{n \in \mathbb{N}} b_n^+ B^n b_{\lambda_1} + \sum_{n \in \mathbb{N}} b_n^- \overline{B}^n \overline{b_{\lambda_1}}$$

and thus the matrix can be expressed in a simple way. After calculation, and denoting $\phi_0 = \sum_{n \in \mathbb{Z}} a_n B^n$, $\phi^+ = \sum_{n \in \mathbb{N}} b_n^+ B^n b_{\lambda_1}$ and $\phi^- = \sum_{n \in \mathbb{N}} b_n^- \overline{B}^n \overline{b_{\lambda_1}}$, we obtain

$$\begin{pmatrix} T_{a z \varphi^+ + \varphi_0 + \bar{a}\bar{z}\varphi^-, K_{z^n}} & T_{b z \varphi^+ + \bar{c}\varphi^-, K_{z^n}} \\ T_{c \varphi^+ + \bar{b}\bar{z}\varphi^-, K_{z^n}} & T_{d \varphi^+ + \varphi_0 + \bar{d}\varphi^+, K_{z^n}} \end{pmatrix}$$

where $a = \overline{b_{\lambda_2}(\lambda_1)}$, $b = <b_{\lambda_1}k_{\lambda_2}, b_{\lambda_2}k_{\lambda_1}> = <k_{\lambda_1}, k_{\lambda_2}>$, $c = <k_{\lambda_1}, k_{\lambda_2}>$, $d = b_{\lambda_1}(\lambda_2)$ and $\varphi_0$, $\varphi^+$ and $\varphi^-$ are associated to $\phi_0$, $\phi^+$ and $\phi^-$. Two questions remain, first ; under the hypothesis that $(a_n)$, $(b_n^+)$ and $(b_n^-)$ are three square summable sequences, can we assure that $a\phi^+ + \phi_0 + \overline{az}\phi^- \in L^2$, second ; can we express eigenvalues or singular values of this matrix.

Moreover we think that it should be possible to show that, if $B$ is any finite Blaschke product then the matrix of $T_{\phi, K_{B^n}}$ has Toeplitz blocks but is not block Toeplitz and it can't be of that form in any basis for $L^2$ containing $\{B^n, n \in \mathbb{Z}\}$, since the diagonal blocks are equal if and only if $\phi \in \overline{\mathrm{Span}\{B^n, n \in \mathbb{Z}\}}^{L^2}$.

From another point of view and using the same idea, we can write $\phi$ in an orthogonal form as

$$\phi = \sum_{n \in \mathbb{Z}} a_n B^n + \sum_{n \in \mathbb{N}} b_n^+ B^n (K_{\lambda_1} - K_{\lambda_2}) + \sum_{n \in \mathbb{N}} b_n^- \overline{B}^n \overline{(K_{\lambda_1} - K_{\lambda_2})}$$

and proposition 11 implies that $\phi \in L^2$ if and only if $(a_n)_{n \in \mathbb{Z}}$, $(b_n^+)_{n \in \mathbb{N}}$ and $(b_n^-)_{n \in \mathbb{N}}$ are three square summable sequences but the matrix is much less simple.

Finally even if the problem of finding a nice basis for $L^2$ which completes $\{B^n, n \in \mathbb{Z}\}$ has many simple answers, the problem is to get a Szegö-type result.